\documentclass[12pt]{article}
\usepackage{xcolor}
\usepackage{amsmath,amsthm,amsfonts,amssymb, mathtools}
\usepackage{needspace}
\usepackage[hidelinks]{hyperref}
\newtheorem{hypA}{Hypothesis}

\numberwithin{equation}{section}

\newtheorem{theorem}{Theorem}[section]
\newtheorem{corollary}[theorem]{Corollary}
\newtheorem{lemma}[theorem]{Lemma}
\newtheorem{remark}[theorem]{Remark}

\def\bC{\mathbb{C}}

\def\bR{\mathbb{R}}

\newcommand{\ud}{\ensuremath{\mathrm{d}}}

\begin{document}
\title{Exact exponential growth for one-dimensional wave equations
with time-dependent potentials}
\author{Raluca M. Balan\footnote{Corresponding author. University of Ottawa, Department of Mathematics
and Statistics, 150 Louis Pasteur Private, Ottawa, Ontario, K1N 6N5, Canada.
E-mail address: rbalan@uottawa.ca.} \and Le Chen\footnote{Auburn University,
Department of Mathematics and Statistics, Office 3147, 305 W. Samford Avenue,
Auburn, Alabama 36849, USA. E-mail address: le.chen@auburn.edu.} \and Jinxin
Wang\footnote{University of Ottawa, Department of Mathematics and Statistics,
150 Louis Pasteur Private, Ottawa, Ontario, K1N 6N5, Canada. E-mail address:
jwang023@uottawa.ca.} }

\date{} \maketitle
\begin{abstract}
  \noindent In this article, we study the asymptotic behaviour of the (mild) solution to the one-dimensional wave equation $v_{tt}=R^2(v_{xx}+Vv)$ with initial condition one and initial velocity zero, as $R \to \infty$. We assume that the potential $V$ is non-negative, bounded, continuously differentiable in space at each time, and uniformly continuous in time, uniformly in space. Our main result shows that the solution $v_R(V;\cdot)$ evaluated at time $t=1$ and location $x=0$ satisfies
  $\lim_{R \to \infty}\frac{1}{R} \log v_R(V;1,0)=\Phi(V)$, for a certain variational constant $\Phi(V)$. Moreover, we provide an upper bound estimate for the same quantity, which holds uniformly for all potentials $V$ with a common bound and a common modulus of continuity in time, for $R$ large enough. These results lie at the core of the developments in the companion article~\cite{balan.chen.ea:26:exact}.
\end{abstract}



\section{Introduction}

In this article, we study the following one-dimensional wave equation with a
time-depend\-ent potential:
\begin{align}
  \label{HAM-v}
  \begin{cases}
    \dfrac{\partial^2 v}{\partial t^2}(t,x)
    = R^2\dfrac{\partial^2 v}{\partial x^2}(t,x)
    +R^2v(t,x)V(t,x),
    \quad t\in[0,1],\ x\in\bR,\\
    v(0,x)=1,\quad
    \dfrac{\partial v}{\partial t}(0,x)=0,
    \quad x\in\bR,
  \end{cases}
\end{align}
where $R>0$ is the wave speed and $V:[0,1]\times \bR \to [0,\infty)$ is a bounded measurable potential. We study the mild solution $v(t,x)\coloneqq v_R(V;t,x)$, which satisfies the integral equation:
\begin{equation}
  \label{duhamel-v}
  \begin{aligned}
    v(t,x)
    &=1+R^2\int_0^t\int_{\bR}
    G_R(t-s,x-y)v(s,y)V(s,y)\,\ud y\,\ud s,
  \end{aligned}
\end{equation}
where $G_R(t,x) \coloneqq (2R)^{-1}\mathbf{1}_{\{|x|<Rt\}}$ for $t>0$ is the
one-dimensional wave kernel. The existence and uniqueness of the solution follow from
the standard Picard iteration (see Section~\ref{S:prelim}).

The simplest example is a constant potential $V\equiv c>0$,
for which the solution of~\eqref{HAM-v} is obtained by averaging
the two exponential solutions in Polyanin and
Nazaikinskii~\cite[Section~6.1.3, p.~567]
{polyanin.nazaikinskii:16:handbook}:
\[
  v_R(c;t,x)=\cosh(Rt\sqrt{c}).
\]

In particular,
\[
  \lim_{R\to\infty}\frac{1}{R}
  \log v_R(c;1,0)=\sqrt{c}.
\]
For a potential varying in both time and space,
such an explicit expression is generally unavailable.
This leads to the following question:%

\begin{center}
\em
What is the exact exponential growth rate of the
solution $v_R(V;1,0)$ as $R\to\infty$?
\end{center}

This question has arisen naturally in the companion
paper~\cite{balan.chen.ea:26:exact}, dedicated to a
stochastic analogue of equation~\eqref{HAM-v} in which
the potential $V$ is replaced by the formal derivative
of a Gaussian noise.

\medskip

We introduce the following hypothesis.

\begin{hypA}\label{hypA}
  The measurable potential $V:[0,1]\times\bR\to[0,\infty)$ satisfies
  \begin{itemize}
    \item[(a)] $V(t,x)\le B$ for all $(t,x)$ and some constant $B>0$ (called a {\em bound});
    \item[(b)] $V(t,\cdot)\in C^1(\bR)$
      for every $t\in[0,1]$;
    \item[(c)] there exists a non-decreasing function $w:[0,1]\to[0,\infty)$ (called a {\em time modulus}) such that $w(0)=0$, $\lim_{h\downarrow0}w(h)=0$, and
      \begin{equation}\label{time-modulus}
        \|V(t,\cdot)-V(s,\cdot)\|_\infty\le w(|t-s|),\qquad \text{for all $s,t\in[0,1]$.}
      \end{equation}
  \end{itemize}
\end{hypA}

For a bounded nonnegative measurable function $U$ on $\bR$, set
\begin{equation}\label{def-lambda}
  \lambda(U) \coloneqq \sup_{\substack{g\in H^1(\bR)\\ \|g\|_2=1}}
  \left\{\int_{\bR}U(x)g(x)^2\,\ud x-\|g'\|_2^2\right\},
\end{equation}
where $H^1(\bR)$ is the real Sobolev space of order one. Suppose $V$ satisfies
Hypothesis~\ref{hypA}. At each time $t$, we apply~\eqref{def-lambda} to the spatial
potential $x \mapsto V(t,x)$. We then define
\begin{equation}\label{def-Phi-V}
  \Phi(V) \coloneqq \int_0^1\sqrt{\lambda(V(t,\cdot))}\,\ud t.
\end{equation}

By Lemma~\ref{lambda-properties}, $\lambda(V(t,\cdot))$ lies in $[0,B]$ and
$|\lambda(V(t,\cdot))-\lambda(V(s,\cdot))|\le w(|t-s|)$. Thus, the integrand
in~\eqref{def-Phi-V} is continuous and the integral is well defined.

We are now ready to state the main result of this article.

\begin{theorem}\label{log-v-th}
  \renewcommand{\theenumi}{(\roman{enumi})}
  \renewcommand{\labelenumi}{\theenumi}
  \begin{enumerate}
    \item\label{fixed-potential-limit} For each fixed $V$ satisfying
    Hypothesis~\ref{hypA},
    \[
      \lim_{R\to\infty}\frac{1}{R}\log v_R(V;1,0)=\Phi(V).
    \]
    \item\label{uniform-upper-bound} Fix $B>0$ and a time modulus $w$ as in
    Hypothesis~\ref{hypA}(c). For every $\varepsilon>0$, there is a threshold
    $R_0=R_0(B,w,\varepsilon)>0$ such that, for every $V$ satisfying
    Hypothesis~\ref{hypA} with bound $B$ and time modulus $w$,
    \[
      \frac{1}{R}\log v_R(V;1,0)\le\Phi(V)+\varepsilon
      \qquad\text{for all }R\ge R_0.
    \]
  \end{enumerate}
\end{theorem}

To prove the upper bound, we divide the time interval into short pieces. On each
piece we hold the potential fixed in time and bound the energy of the wave. For
the lower bound, we measure the size of the wave by a weighted average on a
bounded spatial interval. We choose the weight on each time interval to be a
positive Dirichlet eigenfunction for the potential held fixed there. When we
change weights between intervals, we retain a fixed fraction of the preceding
average. For a fixed partition and spatial interval, this fraction is
independent of $R$, so it has no effect on the logarithmic rate. The wave
equation then lets us bound the value at the origin from below by the last such
average.

The upper bound also explains the role of this result in the study of the
one-dimensional hyperbolic Anderson model~\cite{balan.chen.ea:26:exact}. In that problem the
potential itself changes with $R$. A limit for each fixed potential does not say
how large $R$ must be when a different potential is used at each speed.
Theorem~\ref{log-v-th}\ref{uniform-upper-bound} gives one threshold for the
whole family of potentials with a common bound and time modulus.

\begin{remark}[A potential independent of time]\label{R:time-independent}
  {\rm Suppose $V(t,x)=U(x)$, where $U\in C^1(\bR)$ is bounded and nonnegative.
  Let $u$ solve $u_{tt}=u_{xx}+U(x)u$ with initial displacement one and
  initial velocity zero. Rescaling time gives $v_R(V;t,x)=u(Rt,x)$.
  Since $\lambda(V(s,\cdot))=\lambda(U)$ for every $s$,
  Theorem~\ref{log-v-th}\ref{fixed-potential-limit} gives the long-time rate
  \[
    \lim_{T\to\infty}\frac{1}{T}\log u(T,0)=\sqrt{\lambda(U)}.
  \]

  Separation of variables explains the square root: a spatial mode with
  positive eigenvalue $\lambda$ has a growing solution proportional to
  $e^{t\sqrt{\lambda}}$.

  }
\end{remark}

\begin{remark}[A potential depending only on time]\label{R:time-only}
  {\rm Let $q:[0,1]\to[0,\infty)$ be continuous and set $V(t,x)\coloneqq q(t)$.
  Consider the wave equation with speed $R$, initial displacement one and
  initial velocity zero. Its solution is independent of $x$. Indeed, the
  solution $y_R$ of the ordinary differential equation
  \[
    y_R''(t)=R^2q(t)y_R(t),\qquad y_R(0)=1,\qquad y_R'(0)=0
  \]
  solves the wave equation with these data when extended as a constant
  in $x$. Integrating twice gives
  \[
    y_R(t)=1+R^2\int_0^t(t-s)q(s)y_R(s)\,\ud s,
  \]
  which is~\eqref{duhamel-v} for a solution constant in $x$. By
  Lemma~\ref{L:picard}, $v_R(V;t,x)=y_R(t)$.
  \renewcommand{\theenumi}{(\roman{enumi})}
  \renewcommand{\labelenumi}{\theenumi}
  \begin{enumerate}
    \item \emph{The growth rate.} For a unit test function $g\in H^1(\bR)$,
      the expression in~\eqref{def-lambda}
      is $q(t)-\|g'\|_2^2$. Its supremum is $q(t)$, since dilating $g$ while keeping
      its $L^2$ norm fixed makes $\|g'\|_2$ arbitrarily small. Thus,
      Theorem~\ref{log-v-th}\ref{fixed-potential-limit} and~\eqref{def-Phi-V} give
      \[
        \lim_{R\to\infty}\frac{1}{R}\log y_R(1)
        =\Phi(V)=\int_0^1\sqrt{q(s)}\,\ud s.
      \]
    \item \emph{Olver's approximation.} When $q$ is strictly positive and
      twice continuously differentiable on
      $[0,1]$, Olver's formula gives the same integral. By
      Olver~\cite[Chapter~6, Theorem~2.1]{olver:97:asymptotics}, the ordinary
      differential equation has two solutions that can be normalized to have
      leading terms
      \[
        q(t)^{-1/4}\exp\!\left(\pm R\int_0^t\sqrt{q(s)}\,\ud s\right)
        \qquad (R\to\infty);
      \]
      see Chapter~6, Section~5.1, formula~(5.03). At time one, the coefficient
      of $R$ in the growing exponential is exactly $\Phi(V)$. This agrees with the rate in part (i)
      above for the prescribed initial data. For a potential
      that also varies in space, formula~\eqref{def-Phi-V} replaces $q(s)$
      under the square root by $\lambda(V(s,\cdot))$.
  \end{enumerate}
  }
\end{remark}

\medskip
Wave equations with time-dependent potentials have also been
studied in various settings. Dalang, Mueller and
Tribe~\cite{dalang.mueller.tribe:08:feynman} obtained
Feynman--Kac-type representations for solutions in dimensions $d=1,2,3$.
Colombini et al.~\cite{colombini.petkov.rauch:09:exponential}
proved that for certain nonpositive potentials that are
periodic in time and compactly supported in space, the wave
equation in dimension $d=3$ has solutions with exponentially
growing energy.%

\medskip

For background on wave formula, Duhamel's principle
and energy estimates, see
Evans~\cite[Sections~2.4.1--2.4.3]{evans:10:partial};
for the variational characterization and positivity
of ground states, see Lieb and
Loss~\cite[Theorems~11.5, 11.8 and 12.1]{lieb.loss:01:analysis}.
\medskip

This paper is organized as follows. In Section~\ref{S:prelim}, we construct the solution by Picard iteration, and obtain a preliminary upper bound. In Section~\ref{S:upper}, we prove the uniform upper bound in Theorem~\ref{log-v-th}\ref{uniform-upper-bound}. In Section~\ref{S:lower}, we prove the matching lower bound for each fixed potential. Together the two bounds give the limit in Theorem~\ref{log-v-th}\ref{fixed-potential-limit}.

\section{Preliminaries}
\label{S:prelim}

In this section, we prove that equation~\eqref{HAM-v} has a unique solution, and give a preliminary upper bound for this solution.

\medskip

For bounded measurable $V:[0,1]\times\bR\to\bC$, iteration of
\eqref{duhamel-v} gives a series in ordered space-time variables.
We allow complex values here because the series is applied to a complex-valued
potential in~\cite{balan.chen.ea:26:exact}.
Write $\xi_j\coloneqq (s_j,x_j)$ and $(s_{n+1},x_{n+1})\coloneqq (t,x)$. Define
\begin{equation}\label{def-fnR}
  f_{n,R}(\xi_1,\ldots,\xi_n,t,x)
  \coloneqq\mathbf{1}_{\{0<s_1<\cdots<s_n<t\}}
  \prod_{j=1}^n G_R(s_{j+1}-s_j,x_{j+1}-x_j),
\end{equation}
with $G_R=0$ at nonpositive times. The Picard series is
\begin{equation}
  \label{chaos-v}
  v_R(V;t,x)
  =1+\sum_{n\geq1}R^{2n}
  \int_{([0,1]\times\bR)^n}
  f_{n,R}(\xi_1,\ldots,\xi_n,t,x)
  \prod_{j=1}^n V(\xi_j)\,\ud \xi_1\ldots\,\ud \xi_n.
\end{equation}

\begin{remark}
  \label{support-kernel}
  {\rm Write $\xi_j\coloneqq (s_j,x_j)$, and set $s_{n+1}\coloneqq t$ and $x_{n+1}\coloneqq x$. By
  \eqref{def-fnR} and the support of $G_R$, $f_{n,R}(\xi_1,\ldots,\xi_n,t,x)$
  vanishes unless
  \[
    0<s_1<\cdots<s_n<t,
    \qquad
    |x_{k+1}-x_k|<R(s_{k+1}-s_k),
    \quad k=1,\ldots,n.
  \]
  On this set, $|x-x_j|
  \le \sum_{k=j}^n|x_{k+1}-x_k|
  < R\sum_{k=j}^n(s_{k+1}-s_k)
  =R(t-s_j)\le Rt$ for every $j=1,\ldots,n$.
  Consequently, $f_{n,R}(\cdot,t,x)$ is supported in
  $\bigl([0,t]\times[x-Rt,x+Rt]\bigr)^n$. In particular, $f_{n,R}(\cdot,1,0)$
  is supported in $\bigl([0,1]\times[-R,R]\bigr)^n$. Consequently,
  \begin{equation}
    \label{int-CR}
    v_R(V;1,0)=1+\sum_{n\ge 1} R^{2n}\int_{([0,1]\times[-R,R])^n}
    f_{n,R}(\xi_1,\ldots,\xi_n,1,0)\prod_{j=1}^n V(\xi_j)\,\ud \xi_1\cdots\,\ud \xi_n.
  \end{equation}
  Therefore, $v_R(V;1,0)$
  depends only on the values $V(\xi)$ for $\xi\in[0,1]\times[-R,R]$.
  }
\end{remark}

\begin{lemma}\label{L:picard}
  For $R > 0$ and bounded measurable $V:[0,1]\times\bR \to \bC$,
  equation \eqref{duhamel-v} has a unique bounded solution. It is given by the
  absolutely and uniformly convergent series \eqref{chaos-v}.
  If $V \ge 0$, the solution satisfies: for any $R>0$ and
  $(t,x)\in[0,1]\times\bR$,
  \begin{equation}
    \label{UB-vR}
    1\le v_R(V;t,x)\le\cosh\bigl(Rt\sqrt{\|V\|_\infty}\bigr).
  \end{equation}
\end{lemma}

\begin{proof}
  To see convergence, let $I_n(t)$ be the integral of $f_{n,R}$ over the $n$
  space-time variables, with $I_0(t) \coloneqq 1$. Translation invariance makes this
  integral independent of $x$. Since $\int_{\bR}G_R(s,y)\,\ud y = s$ for $s > 0$,
  integration at the last vertex gives
  \[
    I_{n+1}(t)=\int_0^t(t-s)I_n(s)\,\ud s,\qquad I_n(t)=\frac{t^{2n}}{(2n)!},
  \]
  where the second identity follows by induction. The absolute value of
  the $n$th term is bounded by $(R^2\|V\|_\infty)^n I_n(t)$, so the series
  converges uniformly and substitution in~\eqref{duhamel-v} gives a solution.
  For two bounded solutions the same iteration bounds their difference at
  $(t,x)$ by its supremum norm times $(R^2\|V\|_\infty)^n I_n(t)$ for every
  $n$. Letting $n \to \infty$ proves uniqueness. Summing the majorant by the Taylor series gives
  $\cosh(Rt\sqrt{\|V\|_\infty})$; for nonnegative $V$ all terms are nonnegative.
\end{proof}

\section{Uniform upper bound}\label{S:upper}

In this section, we give the proof of part (ii) of Theorem~\ref{log-v-th}.

\medskip

Throughout this section, all $L^2$ norms and inner products are taken over
$\bR$. We first establish the following properties of $\lambda$.

\begin{lemma}
  \label{lambda-properties}
  For bounded nonnegative spatial potentials $U$ and $\widetilde U$,
  \begin{equation}
    \label{lambda-bounds}
    0\le \lambda(U)\le \|U\|_\infty,
    \qquad
    |\lambda(U)-\lambda(\widetilde U)|
    \le \|U-\widetilde U\|_\infty.
  \end{equation}
\end{lemma}

\begin{proof}
  By \eqref{def-lambda}, we have $\lambda(U)\le \|U\|_\infty,\qquad \lambda(U)-\lambda(\widetilde U) \le \|U-\widetilde U\|_\infty$.
  Interchanging $U$ and $\widetilde U$ gives the second inequality
  in \eqref{lambda-bounds}.

  To prove that $\lambda(U)\ge0$, fix $g\in C_c^\infty(\bR)$
  with $\|g\|_2=1$ and set $g_\ell(x)\coloneqq \ell^{-1/2}g(x/\ell)$.
  Since $\|g_\ell\|_2=1$ and $U\ge0$, we have $\lambda(U)\ge -\|g_\ell'\|_2^2 =-\ell^{-2}\|g'\|_2^2\to0 \text{ as }\ell\to\infty$.
\end{proof}

Suppose that $V$ satisfies Hypothesis~\ref{hypA} with bound $B$ and time modulus $w$.
Fix an integer $m\ge1$. Set $h\coloneqq 1/m$ and $t_i\coloneqq i/m$ for $i=0,\ldots,m$, and
define
\[
  V_i(x)\coloneqq V(t_{i-1},x)+w(h),\qquad
  \lambda_i\coloneqq \lambda(V_i),\qquad i=1,\ldots,m.
\]
Let $\widetilde V(t,x)\coloneqq V_i(x)$ for $t\in[t_{i-1},t_i)$, $i=1,\ldots,m$, and set
$\widetilde V(1,x)\coloneqq V_m(x)$. This potential is piecewise constant in time.

Choose a fixed $\chi\in C_c^\infty(\bR)$ such that
$0\le\chi\le1$, $\chi=1$ on $[-1,1]$, and
$\operatorname{supp}\chi\subset[-2,2]$. Define
\[
  \chi_R(x)\coloneqq \chi\left(\frac{x}{R+1}\right).
\]

Replacing $V$ by $\widetilde V$ and the initial displacement $1$ by $\chi_R$ in
\eqref{HAM-v}, we consider the wave equation
\begin{equation*}
  \begin{cases}
    \dfrac{\partial^2 u}{\partial t^2}(t,x)
    =R^2\dfrac{\partial^2 u}{\partial x^2}(t,x)
    +R^2\widetilde V(t,x)u(t,x),
    & (t,x)\in(t_{i-1},t_i)\times\bR,\quad 1\le i\le m,\\[4pt]
    u(0,x)=\chi_R(x),\quad
    \dfrac{\partial u}{\partial t}(0,x)=0,
    & x\in\bR.
  \end{cases}
\end{equation*}
We study this equation through its Duhamel formula. More precisely, we
seek a bounded function $u$ satisfying, for every $(t,x)\in[0,1]\times\bR$,
\begin{equation}
  \label{duhamel-u}
  u(t,x)=\frac{\chi_R(x+Rt)+\chi_R(x-Rt)}{2}+R^2\int_0^t\int_{\bR}G_R(t-s,x-y)
  \widetilde V(s,y)u(s,y)\,\ud y\,\ud s.
\end{equation}
We next prove that \eqref{duhamel-u} has a unique bounded solution and establish the regularity and support properties needed for the energy estimate on each interval.
\begin{lemma}
  \label{regularity-u}
  For every $R>0$, the integral equation \eqref{duhamel-u}
  has a unique bounded solution $u$. Moreover,
  $u\in C^1([0,1]\times\bR)$, and for every $t\in[0,1]$,
  \[
    \operatorname{supp}u(t,\cdot)
    \subset [-2(R+1)-Rt,\,2(R+1)+Rt].
  \]
  For each $i=1,\ldots,m$, the solution $u$ is $C^2$ on
  $(t_{i-1},t_i)\times\bR$ and satisfies $u_{tt}(t,x)=R^2\bigl(u_{xx}(t,x)+V_i(x)u(t,x)\bigr), \qquad (t,x)\in(t_{i-1},t_i)\times\bR$.
\end{lemma}
\begin{proof}
  Let $u_0(t,x)\coloneqq \frac{\chi_R(x+Rt)+\chi_R(x-Rt)}{2}$.
  Existence and uniqueness follow from the same Picard iteration
  argument as in Lemma~\ref{L:picard}, using $|u_0|\le1$ and
  $\|\widetilde V\|_\infty\le B+w(h)$.
  Since $\chi_R$ is supported in $[-2(R+1),2(R+1)]$,
  the support property follows from the finite propagation
  of the wave kernel.

  By \eqref{duhamel-u} and dominated convergence, $u$ is continuous.
  Since $\widetilde V(s,\cdot)u(s,\cdot)$ is continuous for each
  fixed $s$, Leibniz's rule gives
  \begin{equation}\label{derivatives-u}
    \begin{aligned}
      u_x(t,x)
      &=\partial_x u_0(t,x)
      +\frac R2\int_0^t
      \bigl[(\widetilde V u)(s,x+R(t-s))
      -(\widetilde V u)(s,x-R(t-s))\bigr]\,\ud s,\\
      u_t(t,x)
      &=\partial_t u_0(t,x)
      +\frac{R^2}{2}\int_0^t
      \bigl[(\widetilde V u)(s,x+R(t-s))
      +(\widetilde V u)(s,x-R(t-s))\bigr]\,\ud s.
    \end{aligned}
  \end{equation}
  By dominated convergence, both derivatives are continuous,
  and hence $u\in C^1([0,1]\times\bR)$.

  Since $u$ and $u_x$ have common compact spatial support and
  $V_i\in C^1(\bR)$ for each $i$, the function
  $\partial_x(\widetilde V u)$ is bounded uniformly over the finitely many slabs and continuous
  in space for each fixed time.
  Differentiating the first identity in \eqref{derivatives-u} gives
  \begin{equation}
    \label{uxx}
    u_{xx}(t,x)
    =\partial_x^2u_0(t,x)
    +\frac R2\int_0^t
    \bigl[\partial_x(\widetilde V u)(s,x+R(t-s))
    -\partial_x(\widetilde V u)(s,x-R(t-s))\bigr]\,\ud s.
  \end{equation}
  The mixed derivatives follow by the same rule, with the sum of the two
  boundary values in place of their difference. In $\partial_t u_x$ the two
  terms at $s=t$ cancel, giving $u_{xt}=u_{tx}$. Dominated convergence
  gives continuity of these derivatives and of $u_{xx}$.

  Moreover, for $t\in(t_{i-1},t_i)$, differentiating the
  expression for $u_t$ and using
  $\partial_t^2u_0=R^2\partial_x^2u_0$ gives
  \[
    u_{tt}(t,x)=R^2\bigl(u_{xx}(t,x)+V_i(x)u(t,x)\bigr).
  \]
  The right-hand side is continuous on
  $(t_{i-1},t_i)\times\bR$, which proves the desired
  $C^2$ regularity on each interval.
\end{proof}

\begin{lemma}
  \label{comparison-u-v}
  The solution $u$ of \eqref{duhamel-u} satisfies
  \begin{equation}
    \label{bound-v-by-u}
    1\le v_R(V;1,0)\le u(1,0).
  \end{equation}
\end{lemma}

\begin{proof}
  Iterating \eqref{duhamel-u} at $(1,0)$ gives \eqref{chaos-v} for
  $\widetilde V$, with free term $u_0(1,0)=1$ and a factor $u_0(s_1,x_1)$
  in each term of order $n\ge1$. By Remark~\ref{support-kernel},
  $|x_1|\le R(1-s_1)$ on the kernel support, so $|x_1\pm Rs_1|\le R$.
  Since $\chi_R=1$ on $[-R,R]$, each such factor is one.
  Hence $u(1,0)=v_R(\widetilde V;1,0)$.
  Moreover, by \eqref{time-modulus}, $0\le V(t,x)\le V(t_{i-1},x)+w(h)=\widetilde V(t,x)$ for $t\in[t_{i-1},t_i)$.
  Since $f_{n,R}\ge0$, termwise comparison in
  \eqref{chaos-v} yields $1\le v_R(V;1,0)\le u(1,0)$.
\end{proof}

\begin{proof}[Proof of Theorem~\ref{log-v-th}\ref{uniform-upper-bound}]
  Fix $\eta>0$. For $t\in[t_{i-1},t_i]$, define
  \begin{equation}
    \label{energy}
    E_i(t)
    \coloneqq \frac{1}{R^2}\|u_t(t,\cdot)\|_2^2+\|u_x(t,\cdot)\|_2^2
    +\int_{\bR}(2\lambda_i+\eta-V_i(x))u(t,x)^2\,\ud x.
  \end{equation}
  By Lemma~\ref{regularity-u}, $u(t,\cdot)\in C_c^1(\bR)\subset H^1(\bR)$.
  Applying the definition of $\lambda_i$ to $u(t,\cdot)/\|u(t,\cdot)\|_2$
  when $u(t,\cdot)\ne0$ gives
  $\int_{\bR}V_i(x)u(t,x)^2\,\ud x-\|u_x(t,\cdot)\|_2^2 \le \lambda_i\|u(t,\cdot)\|_2^2$.
  Hence
  \begin{equation}
    \label{energy-lower-bound}
    E_i(t)\ge \frac{1}{R^2}\|u_t(t,\cdot)\|_2^2
    +(\lambda_i+\eta)\|u(t,\cdot)\|_2^2.
  \end{equation}
  In particular, since $\lambda_i\ge0$, we have
  \begin{equation}
    \label{u-energy}
    \|u(t,\cdot)\|_2^2\le \frac{E_i(t)}{\eta}.
  \end{equation}

  By the definition of $E_i$,
  \begin{equation}
    \label{ux-energy}
    \begin{aligned}
      \|u_x(t,\cdot)\|_2^2
      &\le E_i(t)+\int_{\bR}V_i(x)u(t,x)^2\,\ud x\\
      &\le E_i(t)+\|V_i\|_{\infty}\|u(t,\cdot)\|_2^2 \\
      &\le \left(1+\frac{B+w(h)}{\eta}\right)E_i(t).
    \end{aligned}
  \end{equation}
  The last inequality holds because
  $\|V_i\|_\infty\le B+w(h)$.

  For $t\in(t_{i-1},t_i)$, by Lemma~\ref{regularity-u},
  differentiating \eqref{energy} and using
  $u_{tt}=R^2(u_{xx}+V_i u)$, we obtain
  \[
    \begin{aligned}
      E_i'(t)
      &=\frac{2}{R^2}\langle u_t,u_{tt}\rangle
      +2\langle u_x,u_{xt}\rangle
      +2\langle (2\lambda_i+\eta-V_i)u,u_t\rangle\\
      &=2\langle u_t,u_{xx}\rangle
      +2\langle u_x,u_{xt}\rangle
      +2(2\lambda_i+\eta)\langle u,u_t\rangle\\
      &=2(2\lambda_i+\eta)\langle u,u_t\rangle,
    \end{aligned}
  \]
  where the last equality follows by integration by parts,
  since $u_t$ and $u_x$ have compact support and $u_{tx}=u_{xt}$.

  By \eqref{energy-lower-bound},
  \[
    \frac{2\sqrt{\lambda_i+\eta}}{R}\|u\|_2\|u_t\|_2
    \le (\lambda_i+\eta)\|u\|_2^2+\frac{1}{R^2}\|u_t\|_2^2
    \le E_i(t).
  \]
  Consequently,
  \[
    E_i'(t)
    \le R\frac{2\lambda_i+\eta}{\sqrt{\lambda_i+\eta}}E_i(t)
    \le 2R\sqrt{\lambda_i+\eta}\,E_i(t).
  \]
  Since $u,u_x,u_t$ are continuous with a common compact spatial support,
  dominated convergence in \eqref{energy} gives continuity of $E_i$ up to
  both endpoints of its block.
  By Gronwall's inequality,
  \begin{equation}
    \label{energy-growth}
    E_i(t_i) \le
    \exp\!\bigl(2Rh\sqrt{\lambda_i+\eta}\bigr)E_i(t_{i-1}).
  \end{equation}
  We next compare $E_i$ and $E_{i+1}$ at $t_i$, for $1\le i<m$. Since $u$,
  $u_x$, and $u_t$ are continuous at $t_i$, we have
  \begin{equation*}
    E_{i+1}(t_i)-E_i(t_i)=2(\lambda_{i+1}-\lambda_i)\|u(t_i,\cdot)\|_2^2
    -\int_{\bR}(V_{i+1}(x)-V_i(x))u(t_i,x)^2\,\ud x.
  \end{equation*}
  By \eqref{lambda-bounds},
  \begin{equation*}
    |\lambda_{i+1}-\lambda_i|\le \|V_{i+1}-V_i\|_\infty
    =\|V(t_i,\cdot)-V(t_{i-1},\cdot)\|_\infty
    \le w(h).
  \end{equation*}
  By \eqref{u-energy},
  \begin{equation}
    \label{energy-jumps}
    E_{i+1}(t_i)
    \le E_i(t_i)+3w(h)\|u(t_i,\cdot)\|_2^2
    \le \left(1+\frac{3w(h)}{\eta}\right)E_i(t_i).
  \end{equation}
  Combining the estimates \eqref{energy-growth} on the $m$ intervals and \eqref{energy-jumps}
  at the $m-1$ interior times yields
  \[
    E_m(1)
    \le E_1(0)
    \left(1+\frac{3w(h)}{\eta}\right)^{m-1}
    \exp\!\left(2R\sum_{i=1}^m h\sqrt{\lambda_i+\eta}\right).
  \]
  We next calculate $E_1(0)$.
  \begin{equation*}
    \begin{aligned}
      E_1(0)
      =\frac{1}{R^2}\|u_t(0,\cdot)\|_2^2+\|u_x(0,\cdot)\|_2^2
      +\int_{\bR}(2\lambda_{1}+\eta-V_{1}(x))u(0,x)^2\,\ud x.
    \end{aligned}
  \end{equation*}
  At time zero, $u(0,\cdot)=\chi_R$ and $u_t(0,\cdot)=0$. Since
  $0\le\lambda_1\le B+w(h)$,
  \begin{equation*}
    E_1(0)
    \le \frac{\|\chi'\|_2^2}{R+1}
    +(2B+2w(h)+\eta)(R+1)\|\chi\|_2^2 \le C_0(1+R),
  \end{equation*}
  where $C_0$ depends only on $B,w,m,\eta$ and $\chi$.

  At time one, $u(1,\cdot)\in C_c^1(\bR)$, so
  \[
    |u(1,0)|^2
    =2\int_{-\infty}^0 u(1,x)u_x(1,x)\,\ud x
    \le 2\|u(1,\cdot)\|_2\|u_x(1,\cdot)\|_2.
  \]
  By \eqref{u-energy} and \eqref{ux-energy}, applied at $t=1$ with $i=m$, we have
  \[
    \begin{aligned}
      |u(1,0)|^2
      &\le 2\|u(1,\cdot)\|_2\|u_x(1,\cdot)\|_2\\
      &\le 2\left(\frac{E_m(1)}{\eta}\right)^{1/2}
      \left(\left(1+\frac{B+w(h)}{\eta}\right)E_m(1)\right)^{1/2}\\
      &=\frac{2\sqrt{\eta+B+w(h)}}{\eta}\,E_m(1).
    \end{aligned}
  \]
  Using the bound for $E_m(1)$ and then the bound for $E_1(0)$, we obtain
  \[
    |u(1,0)|^2\le \frac{2C_0\sqrt{\eta+B+w(h)}}{\eta}
    \left(1+\frac{3w(h)}{\eta}\right)^{m-1}
    (1+R)\exp\left(2Rh\sum_{i=1}^m\sqrt{\lambda_i+\eta}\right).
  \]
  Taking square roots and using \eqref{bound-v-by-u}, we get
  \[
    v_R(V;1,0)\le u(1,0)
    \le C\sqrt{1+R}\,
    \exp\left(Rh\sum_{i=1}^m\sqrt{\lambda_i+\eta}\right),
  \]
  where $C$ depends only on $B,w,m,\eta$ and $\chi$. In particular, $C$ is
  independent of $R$ and $V$.

  Finally, for $t\in[t_{i-1},t_i]$, we have
  $\|V_i-V(t,\cdot)\|_\infty\le 2w(h)$. Lemma~\ref{lambda-properties} gives
  \[
    \sqrt{\lambda_i+\eta}
    \le \sqrt{\lambda(V(t,\cdot))}+\sqrt{2w(h)+\eta}.
  \]
  Integrating over each interval and summing, we get
  \[
    \sum_{i=1}^m h\sqrt{\lambda_i+\eta}
    \le \int_0^1\sqrt{\lambda(V(t,\cdot))}\,\ud t
    +\sqrt{2w(h)+\eta}
    =\Phi(V)+\sqrt{2w(h)+\eta}.
  \]
  Therefore,
  \[
    \frac1R\log v_R(V;1,0)
    \le \Phi(V)+\sqrt{2w(h)+\eta}
    +\frac{\log\!\bigl(C\sqrt{R+1}\bigr)}{R}.
  \]
  Given $\varepsilon>0$, take $\eta=\varepsilon^2/8$ and choose $m$ large
  enough that $2w(1/m)\le\varepsilon^2/8$. Then
  $\sqrt{2w(h)+\eta}\le\varepsilon/2$. With this choice of $m$ and $\eta$, the
  last term tends to zero uniformly in $V$, since $C$ is independent of $V$.
  Thus there exists $R_0=R_0(B,w,\varepsilon)$ such that, for all $R\ge R_0$,
  \[
    \frac1R\log v_R(V;1,0)\le \Phi(V)+\varepsilon.
  \]
\end{proof}

\begin{remark}[Motivation for the choice of energy]
  \label{R:energy-motivation}
   {\rm The energy defined in \eqref{energy} is a modification
  of the classical conserved energy associated with the wave
  equation $u_{tt}=R^2(u_{xx}+V_i u)$:
  \[
    \mathcal E_i(t)
    =\frac{1}{R^2}\|u_t\|_2^2+\|u_x\|_2^2
    -\int_{\bR}V_i(x)u(t,x)^2\,\ud x.
  \]
  Although $\mathcal E_i(t)$ is conserved on $(t_{i-1},t_i)$,
  it is not necessarily nonnegative or coercive.
  Hence, its conservation alone does not provide
  a bound for $\|u\|_2$.

  The three terms in $(2\lambda_i+\eta)\|u\|_2^2$
  serve different purposes. The first term
  $\lambda_i\|u\|_2^2$ makes the classical energy
  nonnegative by the variational definition of $\lambda_i$.
  The second term $\lambda_i\|u\|_2^2$ is introduced
  to optimize the exponential growth bound.
  Finally, the term $\eta\|u\|_2^2$ ensures
  coercivity uniformly in $\lambda_i$, including
  the case $\lambda_i=0$.
}
\end{remark}

As an immediate consequence, we obtain the upper bound for a fixed potential.

\begin{corollary}
  \label{upper-bound}
  Suppose that $V$ satisfies Hypothesis~\ref{hypA}. Then
  \[
    \limsup_{R\to\infty}\frac{1}{R}\log v_R(V;1,0)\le\Phi(V).
  \]
\end{corollary}
\begin{proof}
  For fixed $V$, take the limsup as $R\to\infty$ in
  Theorem~\ref{log-v-th}\ref{uniform-upper-bound}. Then let $\varepsilon\downarrow0$.
\end{proof}

\section{Lower bound for a fixed potential}\label{S:lower}

In this section, we give the proof of the lower bound in Theorem~\ref{log-v-th}\ref{fixed-potential-limit}.

\medskip

Throughout this section, we assume that $V$ satisfies Hypothesis~\ref{hypA}
and write $v(t,x)\coloneqq v_R(V;t,x)$.
We first prove a comparison lemma.

\begin{lemma}
  \label{cone-comparison}
  For every $R>0$,
  \begin{equation}
    \label{cone-inequality}
    v(s,y)\le v(t,x)
  \end{equation}
  whenever $0\le s\le t\le1$ and $|x-y|\le R(t-s)$.
  In particular, $t\mapsto v(t,x)$ is non-decreasing
  for every $x\in\bR$.
\end{lemma}

\begin{proof}
  By \eqref{chaos-v}, $v\ge1$, and hence $Vv\ge0$.
  If $|x-y|\le R(t-s)$, then for every $0\le r<s$
  and $|z-y|<R(s-r)$, $|z-x|\le |z-y|+|y-x| <R(s-r)+R(t-s)=R(t-r)$.
  Thus the backward cone of $(s,y)$ is contained in that
  of $(t,x)$. By \eqref{duhamel-v} and the nonnegativity
  of $Vv$, we obtain $v(s,y)\le v(t,x)$.
  Taking $y=x$ proves the last assertion.
\end{proof}

Fix an integer $m\ge2$ and $L>0$. Set $h\coloneqq 1/m$, $t_i\coloneqq ih$ for
$i=0,\ldots,m$, and $I\coloneqq (-L,L)$. We reset $m,h,t_i$ here and reserve the last
block $[1-h,1]$ for Lemma~\ref{cone-comparison}, so $i$ runs only to $m-1$.
For $i=1,\ldots,m-1$, define
\[
  \lambda_{i,L}
  \coloneqq \sup_{\substack{g\in H_0^1(I)\\ \|g\|_{L^2(I)}=1}}
  \left\{\int_I V(t_{i-1},x)g(x)^2\,\ud x
  -\int_I |g'(x)|^2\,\ud x\right\}.
\]
Any normalized smooth test function supported in $I$ gives a finite lower
bound, and $V\le B$ gives $\lambda_{i,L}\le B$. Thus $\lambda_{i,L}$ is finite.

\begin{lemma}
  \label{positive-eigenfunction}
  For each $i=1,\ldots,m-1$,
  the supremum defining $\lambda_{i,L}$ is attained by a function
  $\psi_i\in H_0^1(I)\cap C^2(\overline I)$ with
  $\|\psi_i\|_{L^2(I)}=1$, satisfying
  \begin{equation}
    \label{dirichlet-equation}
    \psi_i''(x)+V(t_{i-1},x)\psi_i(x)=\lambda_{i,L}\psi_i(x),
    \qquad x\in I,\qquad \psi_i(-L)=\psi_i(L)=0.
  \end{equation}
  Moreover,
  \[
    \psi_i>0\quad\text{on }I,
    \qquad \psi_i'(-L)>0,
    \qquad \psi_i'(L)<0.
  \]
\end{lemma}

\begin{proof}
  This supremum is attained. Indeed, choose a maximizing sequence
  $(g_n)\subset H_0^1(I)$ with $\|g_n\|_{L^2(I)}=1$ such that
  \[
    \int_I V(t_{i-1},x)g_n(x)^2\,\ud x
    -\|g_n'\|_{L^2(I)}^2
    \ge\lambda_{i,L}-\frac1n.
  \]
  Since $V\le B$, we have
  \[
    \|g_n'\|_{L^2(I)}^2
    \le B-\lambda_{i,L}+\frac1n
    \le B-\lambda_{i,L}+1.
  \]
  Thus $(g_n)$ is bounded in $H_0^1(I)$. Since $H_0^1(I)$ is a
  Hilbert space and the
  embedding $H_0^1(I)\hookrightarrow L^2(I)$ is compact, there exist
  $g\in H_0^1(I)$ and a subsequence, still denoted by $(g_n)$, such that
  \[
    g_n\rightharpoonup g\quad\text{in }H_0^1(I),
    \qquad
    g_n\to g\quad\text{in }L^2(I).
  \]
  In particular, $\|g\|_{L^2(I)}=1$, and
  \[
    \left|\int_I V(t_{i-1},x)
    \bigl(g_n(x)^2-g(x)^2\bigr)\,\ud x\right|
    \le B\|g_n-g\|_{L^2(I)}
    \bigl(\|g_n\|_{L^2(I)}+\|g\|_{L^2(I)}\bigr) \to 0.
  \]
  Moreover, since $g_n'\rightharpoonup g'$ in $L^2(I)$, we have
  \[
    \|g_n'\|_{L^2(I)}^2
    =\|g'\|_{L^2(I)}^2
    +2\int_I g'(x)\bigl(g_n'(x)-g'(x)\bigr)\,\ud x
    +\|g_n'-g'\|_{L^2(I)}^2.
  \]
  The integral tends to zero by weak convergence. We obtain
  \[
    \|g'\|_{L^2(I)}^2
    \le\liminf_{n\to\infty}\|g_n'\|_{L^2(I)}^2.
  \]
  Consequently,
  \begin{equation*}
    \lambda_{i,L}
    =\lim_{n\to\infty}
    \left\{\int_I V(t_{i-1},x)g_n(x)^2\,\ud x
    -\|g_n'\|_{L^2(I)}^2\right\}
    \le\int_I V(t_{i-1},x)g(x)^2\,\ud x
    -\|g'\|_{L^2(I)}^2
    \le\lambda_{i,L}.
  \end{equation*}

  The supremum is attained.
  Replacing a maximizer by its absolute value, we obtain a
  nonnegative maximizer $\psi_i$.

  The Euler--Lagrange equation gives
  \[
    \int_I V(t_{i-1},x)\psi_i(x)\varphi(x)\,\ud x
    -\int_I\psi_i'(x)\varphi'(x)\,\ud x
    =\lambda_{i,L}\int_I\psi_i(x)\varphi(x)\,\ud x
  \]
  for every $\varphi\in H_0^1(I)$.
  Since $\psi_i\in H_0^1(I)$,
  it has a continuous representative on $\overline I$. The weak derivative of
  $\psi_i'$ equals the continuous function $(\lambda_{i,L}-V(t_{i-1},\cdot))\psi_i$.
  Integration gives $\psi_i\in C^2(\overline I)$, so the equation holds classically.

  An interior zero of the nonnegative $\psi_i$ would be a minimum, so its derivative would
  vanish there. Uniqueness for \eqref{dirichlet-equation} with these
  initial values would give $\psi_i\equiv0$. Thus $\psi_i>0$ on $I$.
  By the same argument, its derivative cannot vanish at either
  endpoint. Since $\psi_i>0$ on $I$ and $\psi_i(\pm L)=0$, we have
  $\psi_i'(-L)\ge0\ge\psi_i'(L)$. Hence $\psi_i'(-L)>0$ and $\psi_i'(L)<0$.
\end{proof}

\begin{proof}[Proof of Theorem~\ref{log-v-th}\ref{fixed-potential-limit}]
  \setlength{\emergencystretch}{1em}
  Let $\psi_i$ be the positive eigenfunction
  given by Lemma~\ref{positive-eigenfunction}.
  Rescale $\psi_i$ so that
  \[
    \int_I\psi_i(x)\,\ud x=1.
  \]
  This rescaling preserves \eqref{dirichlet-equation} and the signs of
  the endpoint derivatives. Set
  \[
    q_i(t)\coloneqq \int_I v(t,x)\psi_i(x)\,\ud x.
  \]
  Since $\psi_i\ge0$, Lemma~\ref{cone-comparison} implies that $q_i$ is
  non-decreasing. Moreover, $q_i(0)=1$.

  We first derive a differential inequality for $q_i$ on
  $[t_{i-1},t_i]$.
  By Hypothesis~\ref{hypA}(b),(c),
  \[
    |V(t,x)-V(s,y)|\le w(|t-s|)+|V(s,x)-V(s,y)|,
  \]
  so $V$ is jointly continuous.

  Apply the first-derivative argument of Lemma~\ref{regularity-u} with
  $u_0=1$ and $Vv$ in place of $\widetilde Vu$. The integration domains
  are uniformly bounded on compact sets of $(t,x)$, and $Vv$ is bounded.
  Thus
  $v\in C^1([0,1]\times\bR)$, with one-sided time derivatives
  at the endpoints, and
  \begin{equation}\label{derivatives-v}
    \begin{aligned}
      v_t(t,x)
      &=\frac{R^2}{2}\int_0^t
      \bigl[(Vv)(s,x+R(t-s))
      +(Vv)(s,x-R(t-s))\bigr]\,\ud s,\\
      v_x(t,x)
      &=\frac R2\int_0^t
      \bigl[(Vv)(s,x+R(t-s))
      -(Vv)(s,x-R(t-s))\bigr]\,\ud s.
    \end{aligned}
  \end{equation}

  Since $v\in C^1$, \eqref{derivatives-v} gives $q_i'(t)=\int_I v_t(t,x)\psi_i(x)\,\ud x$, that is,
  \[
    q_i'(t)=\frac{R^2}{2}\int_0^t\int_I
    \big[(Vv)(s,x+R(t-s))+(Vv)(s,x-R(t-s))\big]\psi_i(x)\,\ud x\,\ud s.
  \]
  By a change of variables and Leibniz's rule, using
  $\psi_i(\pm L)=0$, we obtain
  \[
    \frac{\ud}{\ud a}\int_I(Vv)(s,x+a)\psi_i(x)\,\ud x
    =-\int_I(Vv)(s,x+a)\psi_i'(x)\,\ud x.
  \]
  Differentiating in $t$, the upper limit $s=t$ contributes
  $R^2\int_I(Vv)(t,x)\psi_i(x)\,\ud x$. Applying the shift identity with
  $a=\pm R(t-s)$ to the remaining terms and using \eqref{derivatives-v} gives
  \[
    q_i''(t)=R^2\int_I
    \bigl(V(t,x)v(t,x)\psi_i(x)
    -v_x(t,x)\psi_i'(x)\bigr)\,\ud x.
  \]
  The right-hand side is continuous in $t$, so
  $q_i\in C^2([0,1])$.

  Integrating by parts in $x$, we obtain
  \begin{equation*}
    \frac{q_i''(t)}{R^2}  =\int_I v(t,x)
    \bigl(\psi_i''(x)+V(t,x)\psi_i(x)\bigr)\,\ud x +v(t,-L)\psi_i'(-L)-v(t,L)\psi_i'(L).
  \end{equation*}
  By Lemma~\ref{positive-eigenfunction}, both boundary terms are
  nonnegative. For $t\in[t_{i-1},t_i]$, adding and subtracting
  $V(t_{i-1},x)\psi_i(x)$ and using \eqref{dirichlet-equation}, we obtain
  \[
    \begin{aligned}
      \frac{q_i''(t)}{R^2}
      &\ge\int_I v(t,x)
      \bigl(\psi_i''(x)+V(t,x)\psi_i(x)\bigr)\,\ud x\\
      &=\int_I v(t,x)
      \bigl(\psi_i''(x)+V(t_{i-1},x)\psi_i(x)\bigr)\,\ud x+\int_I v(t,x)
      \bigl(V(t,x)-V(t_{i-1},x)\bigr)\psi_i(x)\,\ud x\\
      &=\lambda_{i,L}q_i(t)
      +\int_I v(t,x)
      \bigl(V(t,x)-V(t_{i-1},x)\bigr)\psi_i(x)\,\ud x.
    \end{aligned}
  \]
  By \eqref{time-modulus},
  $V(t,x)-V(t_{i-1},x)\ge-w(h)$.
  Since $v$ and $\psi_i$ are nonnegative, it follows that
  \begin{equation}
    \label{inequality-q}
    q_i''(t)\ge R^2\bigl(\lambda_{i,L}-w(h)\bigr)q_i(t).
  \end{equation}

  We next use this inequality to estimate the growth of $q_i$.
  Set
  \[
    \mu_i\coloneqq \sqrt{(\lambda_{i,L}-w(h))_+}.
  \]

  If $\mu_i>0$, then for $t_{i-1}\le a\le b\le t_i$,
  monotonicity gives $q_i'(a)\ge0$ and
  $q_i(r)\ge q_i(a)$ for all $r\in[a,b]$.
  Using \eqref{inequality-q}, we obtain
  \[
    \begin{aligned}
      q_i(b)
      &=q_i(a)+(b-a)q_i'(a)
      +\int_a^b(b-r)q_i''(r)\,\ud r\\
      &\ge q_i(a)+R^2\mu_i^2\int_a^b(b-r)q_i(r)\,\ud r\\
      &\ge q_i(a)+R^2\mu_i^2q_i(a)\int_a^b(b-r)\,\ud r\\
      &=q_i(a)\left(1+\frac{(R\mu_i(b-a))^2}{2!}\right).
    \end{aligned}
  \]
  We claim that for every $n\ge0$ and $t_{i-1}\le a\le r\le t_i$,
  \[
    q_i(r)\ge q_i(a)\sum_{k=0}^{n}\frac{(R\mu_i(r-a))^{2k}}{(2k)!}.
  \]
  The case $n=0$ follows from monotonicity. Suppose the claim holds for $n$.
  Insert this lower bound for $q_i(r)$ into the first inequality above and use
  \[
    \int_a^b (b-r)\frac{(r-a)^{2k}}{(2k)!}\,\ud r
    =\frac{(b-a)^{2k+2}}{(2k+2)!}.
  \]
  This gives the claim for $n+1$ at every $b\in[a,t_i]$.
  Letting $n\to\infty$ yields
  \[
    q_i(b)\ge q_i(a)\sum_{n=0}^{\infty}
    \frac{(R\mu_i(b-a))^{2n}}{(2n)!}
    =q_i(a)\cosh(R\mu_i(b-a)).
  \]
  If $\mu_i=0$, monotonicity gives $q_i(b)\ge q_i(a)$.
  Thus, in both cases,
  we obtain
  \begin{equation}
    \label{growth-q}
    q_i(t_i)\ge\frac12 e^{Rh\mu_i}q_i(t_{i-1}).
  \end{equation}

  To combine these estimates, we compare $q_i$ and $q_{i+1}$ at
  $t_i$, for $1\le i\le m-2$.
  By Lemma~\ref{positive-eigenfunction}, both $\psi_i$ and
  $\psi_{i+1}$ are positive on $I$ and vanish at the endpoints,
  where their derivatives are nonzero. By L'H\^{o}pital's rule,
  \[
    \lim_{x\downarrow-L}\frac{\psi_{i+1}(x)}{\psi_i(x)}
    =\frac{\psi_{i+1}'(-L)}{\psi_i'(-L)}>0,
    \qquad
    \lim_{x\uparrow L}\frac{\psi_{i+1}(x)}{\psi_i(x)}
    =\frac{\psi_{i+1}'(L)}{\psi_i'(L)}>0.
  \]
  Thus $\psi_{i+1}/\psi_i$ extends to a positive continuous function
  on $\overline I=[-L,L]$, so
  \[
    c_i\coloneqq \inf_{x\in I}\frac{\psi_{i+1}(x)}{\psi_i(x)}>0.
  \]
  We obtain
  \[
    q_{i+1}(t_i)
    =\int_I v(t_i,x)\psi_{i+1}(x)\,\ud x
    \ge c_i\int_I v(t_i,x)\psi_i(x)\,\ud x
    =c_iq_i(t_i).
  \]
  Since $q_1(0)=1$, combining this estimate with
  \eqref{growth-q} on the first $m-1$ intervals gives
  \[
    q_{m-1}(1-h)
    \ge 2^{-(m-1)}
    \left(\prod_{i=1}^{m-2}c_i\right)
    \exp\left(Rh\sum_{i=1}^{m-1}\mu_i\right).
  \]

  We now pass from $q_{m-1}(1-h)$ to $v(1,0)$.
  For $Rh\ge L$ and $y\in I$, Lemma~\ref{cone-comparison}
  gives $v(1-h,y)\le v(1,0)$. Multiplying by $\psi_{m-1}(y)$
  and integrating over $I$, we obtain
  \[
    q_{m-1}(1-h)
    \le v(1,0)\int_I\psi_{m-1}(y)\,\ud y
    =v(1,0).
  \]
  For fixed $m$, $L$, and $V$, the constants $c_i$ are independent
  of $R$. Taking logarithms, dividing by $R$, and letting
  $R\to\infty$, we get
  \[
    \liminf_{R\to\infty}\frac1R\log v_R(V;1,0)
    \ge h\sum_{i=1}^{m-1}
    \sqrt{(\lambda_{i,L}-w(h))_+}.
  \]

  It remains to let $L\to\infty$ and then $m\to\infty$.
  We first show that, for each $i=1,\ldots,m-1$,
  \[
    \lambda_{i,L}\uparrow\lambda(V(t_{i-1},\cdot))
    \qquad\text{as }L\to\infty.
  \]

  By extension by zero, $\lambda_{i,L}$ is non-decreasing in $L$
  and bounded above by $\lambda(V(t_{i-1},\cdot))$.
  Conversely, by the density of $C_c^\infty(\bR)$ in $H^1(\bR)$
  and the continuity of the variational functional, the supremum
  in \eqref{def-lambda} can be restricted to normalized functions
  in $C_c^\infty(\bR)$. Since every such function belongs to
  $H_0^1(-L,L)$ for sufficiently large $L$, the desired
  convergence follows.

  Letting $L\to\infty$, we obtain
  \[
    \liminf_{R\to\infty}\frac1R\log v_R(V;1,0)
    \ge h\sum_{i=1}^{m-1}
    \sqrt{\bigl(\lambda(V(t_{i-1},\cdot))-w(h)\bigr)_+}.
  \]

  Finally, for $t\in[t_{i-1},t_i]$, by \eqref{time-modulus}, Lemma~\ref{lambda-properties}, and the inequality
  $\sqrt{(a-b)_+}\ge\sqrt a-\sqrt b$ for $a,b\ge0$, we have
  \[
    \sqrt{\bigl(\lambda(V(t_{i-1},\cdot))-w(h)\bigr)_+}
    \ge\sqrt{\lambda(V(t,\cdot))}-\sqrt{2w(h)}.
  \]
  Integrating over the first $m-1$ intervals and summing, we get
  \[
    \begin{aligned}
      h\sum_{i=1}^{m-1}
      \sqrt{\bigl(\lambda(V(t_{i-1},\cdot))-w(h)\bigr)_+}
      &\ge\int_0^{1-h}\sqrt{\lambda(V(t,\cdot))}\,\ud t
      -(1-h)\sqrt{2w(h)}\\
      &\ge\Phi(V)-h\sqrt B-(1-h)\sqrt{2w(h)},
    \end{aligned}
  \]
  where the last inequality follows from
  $\lambda(V(t,\cdot))\le B$ on $[1-h,1]$.
  Therefore,
  \[
    \liminf_{R\to\infty}\frac1R\log v_R(V;1,0)
    \ge\Phi(V)-(1-h)\sqrt{2w(h)}-h\sqrt B.
  \]
  Letting $m\to\infty$, so that $h\downarrow0$ and $w(h)\to0$, gives
  the lower bound. Combining it with Corollary~\ref{upper-bound} gives the limit in
  Theorem~\ref{log-v-th}\ref{fixed-potential-limit}.
\end{proof}

\section*{Acknowledgments}
Jinxin Wang thanks the Department of Mathematics and Statistics at Auburn
University for its hospitality during his research visit. Le Chen thanks OpenAI for
providing access to its AI tools through the ChatGPT for Academic Researchers
program.

\paragraph{Funding.}
Raluca M.~Balan was supported by a grant from the Natural Sciences and
Engineering Research Council of Canada. Le Chen was partially supported by
the National Science Foundation through CAREER award \mbox{DMS-2443823} and by
the Simons Foundation through award \mbox{MPS-TSM-00959981}. A Mitacs
Globalink Research Award Abroad (IT52612) supported Jinxin Wang's visit to
Auburn University.

\paragraph{Statement on AI use.}
The proof strategy and a
first version of the arguments were obtained by the authors with the help
of an AI system. The validity of all proofs was checked by the authors, who take full
responsibility for the content. Arriving at the argument required many exchanges, in which the authors examined the proposed steps, and asked follow-up questions. Even after the main idea emerged, the argument was difficult to follow in the form that was provided. The authors spent a substantial amount of time and effort reconstructing the proof, checking its details, and deciding how to explain it. The final presentation is the illustration of that work.

\end{document}